\documentclass[letterpaper, 10 pt, conference]{ieeeconf} 

\IEEEoverridecommandlockouts
\usepackage{graphicx}
\usepackage{amsmath}
\usepackage{amssymb}
\usepackage{amsfonts}
\usepackage{latexsym}
\usepackage{epsfig}
\usepackage{theorem}
\usepackage{shadow}
\usepackage{float}
\usepackage{color}
\usepackage{enumerate}
\usepackage{psfrag}
\usepackage{epstopdf}

\newtheorem{lem}{\bf{Lemma}}

\newtheorem{rem}{Remark}
\newtheorem{defi}{Definition}
\newtheorem{hypo}{\bf{Assumption}}

\begin{document}
\title{\bf Riccati observers for velocity-aided attitude estimation of accelerated vehicles using coupled velocity measurements}

\author{Minh-Duc Hua, Tarek Hamel, Claude Samson
\thanks{M.-D. Hua and T. Hamel are with I3S, Universit\'e C\^ote d'Azur, CNRS, Sophia Antipolis, France, {\tt\footnotesize $hua(thamel)@i3s.unice.fr$}.}
\thanks{C. Samson is with INRIA and I3S UCA-CNRS, Sophia Antipolis, France, {\tt\footnotesize $claude.samson@inria.fr$}, {\tt\footnotesize $csamson@i3s.unice.fr$}.}
}

\maketitle

\pagestyle{empty}
\thispagestyle{empty}

\begin{abstract}
Motivated by drone autonomous navigation applications we address a novel problem of velocity-aided attitude estimation by combining two linear velocity components measured in a body-fixed frame and a linear velocity component measured in an inertial frame with the measurements of an Inertial Measurement Unit (IMU). The main contributions of the present paper are the design of Riccati nonlinear observers, which may be viewed as deterministic versions of an Extended Kalman filter (EKF), and an analysis of observability conditions under which local exponential stability of the observer is achieved. Reported simulation results further indicate that the observers' domain of convergence is large.
\end{abstract}

\section{Introduction}\label{sec:introduction}

The development of reliable attitude (i.e. orientation) estimators is a key requirement for efficient automatic control of drones. Most existing attitude observers exploit the measurements of an IMU under the assumption of weak accelerations of the vehicle to justify the direct use of accelerometer measurements for the estimation of the gravity direction in a body-fixed frame \cite{hua14,mhp08,ms10cep, batista2012}. The violation of this assumption, when the vehicle undergoes sustained accelerations, jeopardizes the accuracy of the attitude estimate (cf. \cite{hua10cep}). To overcome this problem velocity-aided attitude observers have been developed in the literature by fusing IMU measurements with the vehicle's linear velocity measurements done either in an inertial frame \cite{grip2012nonlinear,hua10cep,ms08IFAC,rt11CDC, khosravian16} or in a body-fixed frame~\cite{bonnabelITAC08,troniICRA13,hua16,Allibert16}. The present paper addresses a new problem of velocity-aided attitude estimation where the vehicle's linear velocity is measured partly in a body-fixed frame and partly in an inertial frame. A motivating application of this work is related to quadrotor UAV navigation in situations where linear velocity's components along two body axes orthogonal to the thrust direction and expressed in a body-fixed frame can be derived from accelerometer measurements combined with an aerodynamic linear drag model \cite{Martin10ICRA,mahony2012} and where the linear velocity's vertical component expressed in an inertial frame can be obtained from barometer measurements. The important nonlinearities resulting from the use of such measurements render the design of an attitude observer significantly more complex than when all the linear velocity's components are measured in a single frame, either inertial or body-fixed. They also exclude the possibility of proving semi-global, or almost-global, stability results similar to these derived in \cite{hua10cep} and \cite{hua16} in the simpler case of complete linear velocity measurements in a single frame.

The design of the observers proposed in this paper are adapted from a recent deterministic Riccati observer design framework \cite{HamSamTAC16} that relies on the solutions to the Continuous Riccati Equation (CRE) and encompasses EKF solutions. Accordingly, good conditioning of the solutions to the CRE and, subsequently, exponential stability of the obtained  observers rely on conditions of {\em uniform observability} whose satisfaction calls for a specific analysis. Since only local stability is demonstrated simulation results are useful to get complementary indications about the performance and the size of the basin of attraction of these observers.

The paper is organised as follows. Notation, system equations, and the measurements involved in the observer design are specified in Section  \ref{sec:preliminary}. In the same section some basic definitions and conditions about system observability are recalled, together with elements of the deterministic Riccati observer design framework proposed in \cite{HamSamTAC16}. In Section \ref{sec:observerdesign} the observers expressions are specified, and an analysis of associated observability conditions is carried out in Section \ref{sec:observability_analysis}. Simulation results illustrating the performance of the observers and showing that their domain of convergence can be large are reported in Section \ref{sec:simulation}. A short concluding section follows.

\vspace{-0.1cm}
\section{Preliminary material}\label{sec:preliminary}
\subsection{Notation}

\noindent$\bullet\,$ $\{e_1,e_2,e_3\}$ denotes the canonical basis of $\mathbb{R}^3$ and $[\cdot]_\times$ denotes the skew-symmetric matrix associated with the cross product, i.e., $[u]_\times v = u \times v, \forall u, v \in \mathbb{R}^3$. The identity matrix of $\mathbb{R}^{n\times n}$  is denoted as $I_n$
and $\pi_x \triangleq I_3 - x x^\top$, $\forall x \in S^2$ (the unit 2-sphere), is the projection operator onto the plane orthogonal to $x$. Note that $\pi_x = -[x]_\times^2$, $\forall x \in S^2$. \\
$\bullet\,$ $\{\mathcal{I}\} = \{O; \vec \imath_0, \vec\jmath_0,\vec k_0\}$ denotes an inertial frame attached to the earth, typically chosen as the north-east-down frame, and $\{\mathcal B\} = \{G; \vec \imath, \vec\jmath,\vec k\}$ is a body-fixed frame whose origin is the vehicle's center of mass $G$.\\
$\bullet\,$ The vehicle's attitude is represented by a rotation matrix $R \in \text{SO(3)}$ of the frame $\{\mathcal B\}$ relative to  $\{\mathcal I\}$. The column vectors of $R$ correspond to the vectors of coordinates of $\vec \imath,\vec\jmath, \vec k$ expressed in the basis of $\{\mathcal I\}$. The element at the intersection of the $i^{th}$ row and $j^{th}$ column of $R$ is denoted as $R_{i,j}$, with $i,j \in \{1,2,3\}$.\\
$\bullet\,$ $V \in \mathbb{R}^3$ and $\Omega \in \mathbb{R}^3$ are the vectors of coordinates of the vehicle's linear and angular velocities expressed in $\{\mathcal B\}$. The linear velocity expressed in $\{\mathcal I\}$ is denoted as $v \in \mathbb{R}^3$ so that $v = R V$.

\subsection{System equations and measurements}
The vehicle's attitude satisfies the differential equation \vspace{-0.1cm}
\begin{equation}\label{attitudeSystem}
\dot R = R [\Omega]_\times \vspace{-0.1cm}
\end{equation}
We assume that the vehicle is equipped with an IMU consisting of a 3-axis gyrometer that measures the angular velocity $\Omega \in \mathbb{R}^3$ and of a 3-axis accelerometer that measures the {\it specific acceleration} $a_{\mathcal B} \in \mathbb{R}^3$, expressed in $\{\mathcal B\}$.
Using the flat non-rotating Earth assumption, we have \cite{bonnabelITAC08} \vspace{-0.1cm}
\begin{equation}\label{dotv}
\dot V =  -[\Omega]_\times V + a_{\mathcal B} + g R^\top e_3 \vspace{-0.1cm}
\end{equation}
where $g$ is the gravity constant.
A 3-axis magnetometer is also integrated in many IMUs to measure of the normalized Earth's magnetic field vector expressed in $\{\mathcal B\}$. Let $m_{\mathcal I} = [m_1,m_2,m_3]^\top \in S^2$ denote the known normalised Earth's magnetic field vector expressed in $\{\mathcal I\}$. The vectors $m_{\mathcal I}$ and $e_3$ are usually assumed to be non-collinear so that $R$ can be estimated from the observation (measurements) in the body-fixed frame of the gravity vector and of the Earth's magnetic field vector. The magnetometer thus measures $m_{\mathcal B} = R^\top m_{\mathcal I}$.
We further assume that the vehicle is equipped with sensory devices that provide measurements of the two first components of $V$ and the third component of $v$, i.e., $V_1$, $V_2$ and $v_3$. A possible combination of sensors providing such measurements in the case of a flying drone was evoked in the introduction. To summarize, we assume that the available measurements are $V_1$ and $V_2$ (that may be provided by an onboard accelerometer), $v_3$ (that may be provided by an onboard barometer), and $m_{\mathcal B}$ (that is provided by an onboard magnetometer).

\vspace{-0.1cm}
\subsection{Recalls of observability definitions and conditions}
The following definitions and conditions are classical and just recalled here for the sake of completeness.
Consider a linear time-varying (LTV) system \vspace{-0.cm}
\begin{equation}\label{LTVgen}
\left\{\!\!
\begin{array}{ll}
\dot x &\!\!\!\! =A(t) x + B(t) u \\
y &\!\!\!\!= C(t)x
\end{array}
\right.\vspace{-0.cm}
\end{equation}
with $x \in \mathbb{R}^n$ the system state vector, $u \in \mathbb{R}^s$ the system input vector, and $y \in \mathbb{R}^m$ the system output vector.

\vspace{-0.3cm}
\begin{defi} (instantaneous observability) System \eqref{LTVgen} is instantaneously observable if $\forall t$, $x(t)$ can be calculated from the input $u(t)$, the output $y(t)$, and the time-derivatives $u^{(k)}(t), y^{(k)}(t), k\in \mathbb{N}$.
\end{defi}

\vspace{-0.5cm}
\begin{lem}\label{lemObservabilityInstantaneous}
(see \cite{brockett1972}) Define the observation space at $t$ as the space generated by \vspace{-0.2cm}
\[
{\mathcal O}(t) \triangleq \begin{bmatrix} N_0(t) \\ N_1(t) \\ \vdots \end{bmatrix} \vspace{-0.1cm}
\]
with $N_0 = C$, $N_{k} = N_{k-1} A + \dot N_{k-1}, k={1, \cdots}$. Then, System \eqref{LTVgen} is instantaneously observable if $\mathrm{rank}({\mathcal O}) = n$.
\end{lem}

\vspace{-0.5cm}
\begin{defi} (uniform observability) System \eqref{LTVgen} is uniformly observable if there exists $\delta>0$,  $\mu>0$ such that $\forall t\geq 0$ \vspace{-0.1cm}
\begin{equation}\label{grammian}
W(t,t+\delta) \triangleq \frac{1}{\delta} \int_t^{t+\delta} \Phi^\top(t,s) C^\top(s) C(s) \Phi(t,s) ds \vspace{-0.1cm}
\end{equation}
with $\Phi(t,s)$ the transition matrix associated with $A(t)$, i.e. such that $\frac{d}{dt} \Phi(t,s) = A(t) \Phi(t,s)$ with $\Phi(t,t) = I_n$.
\end{defi}
\vspace{-0.1cm}
$W(t,t+\delta)$ is called the observability Gramian of System \eqref{LTVgen}. When \eqref{grammian} is satisfied one also says that the pair $(A(t),C(t))$ is uniformly observable. The following useful condition, derived in \cite{scandaroli2013}, points out a sufficient condition for uniform observability.

\vspace{-0.2cm}
\begin{lem}\label{lemObservabilityUniform}
(see \cite{scandaroli2013}) If there exists a matrix-valued function $M(\cdot)$ of dimension $(p\times n)$ $(p\geq 1)$ composed of row vectors of $N_0(\cdot)$, $N_1(\cdot)$, $\cdots$, such that for some positive numbers $\bar\delta, \bar\mu$ and $\forall t\geq 0$
\begin{equation}\label{condUO}
\frac{1}{\bar\delta} \int_{t}^{t+\bar\delta} M^\top(s) M(s) ds \geq \bar\mu I_n
\end{equation}
then the observability Gramian of System \eqref{LTVgen} satisfies condition \eqref{grammian}.
\end{lem}

\vspace{-0.4cm}
\begin{rem}\label{rem1}
It is noticeable that System \eqref{LTVgen} can be uniformly observable but not instantaneously observable. Instantaneous observability of System \eqref{LTVgen} does not either imply uniform observability. For instance, the matrix $M$ involved in Lemma \ref{lemObservabilityUniform} may be always of full rank with the determinant of $M^\top M$ converging to zero (while remaining positive) and such that condition \eqref{condUO} is not satisfied. However, if there exists $\mu>0$ such that  $M^\top(t) M(t) \geq \mu I_n ,\,\, \forall t\geq 0$
then the system is both instantaneously observable and uniformly observable.
\end{rem}

\vspace{-0.2cm}
\subsection{Recalls of a Riccati observer design framework} \label{subsec:Riccati}
The proposed observer design is adapted from the deterministic observer design framework reported in \cite{HamSamTAC16}. Consider the nonlinear system \vspace{-0.cm}
\begin{equation}\label{LTV}
\left\{\!\!
\begin{array}{ll}
\dot {x} &\!\!\!\!= A(x_1, t) x + u \\
y &\!\!\!\!= C(x, t)x
\end{array}\right. \vspace{-0.cm}
\end{equation}
with $x = [x_1^\top, x_2^\top]^\top$, $x_1 \in {\cal{B}}^n_1$ (the closed unit ball in $\mathbb{R}^n$), $x_2 \in \mathbb{R}^n$, $y \in \mathbb{R}^m$, $A(x_1, t)$ a continuous matrix-valued function uniformly bounded w.r.t. $t$ and uniformly continuous w.r.t. $x_1$ of the form \vspace{-0.1cm}
\[
A(x_1, t) = \begin{bmatrix}
A_{1,1}(t) & 0 \\ A_{2,1}(x_1,t) & A_{2,2}(t)
\end{bmatrix} \vspace{-0.1cm}
\]
and $C(x, t)$ a continuous matrix-valued function uniformly bounded w.r.t. $t$ and uniformly continuous w.r.t. $x$.
Apply the input \vspace{-0.cm}
\begin{equation}\label{U}
u = -k(t) P C^\top Q y \vspace{-0.cm}
\end{equation}
with $0.5 \leq k(t) \leq k_{max} <\infty$ and $P\in \mathbb{R}^{2n\times 2n}$ a symmetric positive definite matrix solution to the following CRE \vspace{-0.cm}
\begin{equation}\label{CRE}
\dot P = AP + P A^\top -P C^\top Q(t) C P + S(t) \vspace{-0.cm}
\end{equation}
with $P(0) \in \mathbb{R}^{2n\times 2n}$ a symmetric positive definite matrix, $Q(t)\in \mathbb{R}^{m\times m}$ bounded continuous symmetric positive semidefinite, and $S(t)\in \mathbb{R}^{2n\times 2n}$ bounded continuous symmetric positive definite. Then, from Theorem 3.1 and Corollary 3.2 in \cite{HamSamTAC16}, $x = 0$ is locally exponentially stable when $Q(t)$ and $S(t)$ are both larger than some positive matrix and the pair $(A^\star(t), C^\star(t))$, with \vspace{-0.1cm}
\begin{equation}\label{ACstar}
A^\star(t) \triangleq A(0,t),\,\, C^\star(t) \triangleq C(0,t) \vspace{-0.1cm}
\end{equation}
is uniformly observable.

\vspace{-0.1cm}
\section{Observer design}\label{sec:observerdesign}
Let $\hat R \in \mathrm{SO(3)}$ and $\hat V \in  \mathbb{R}^3$ denote the estimates of $R$ and $V$, respectively. The proposed observers are of the form \vspace{-0.cm}
\begin{equation}\label{observerGen}
\left\{
\begin{array}{ll}
&\!\!\!\!\!\!\!\! \dot{\hat R} = \hat R [\Omega - \sigma_R]_\times \\
&\!\!\!\!\!\!\!\! \dot{\hat V} = -[\Omega]_\times\hat V + a_{\mathcal B} + g \hat R^\top e_3 - \sigma_V
\end{array}\right. \vspace{-0.cm}
\end{equation}
where $\sigma_R, \sigma_V \in \mathbb{R}^3$ are innovation terms specified thereafter.
Defining the observer errors \vspace{-0.cm}
\[
\tilde R \triangleq R \hat R^\top, \,\, \bar R \triangleq \hat R^\top R, \,\, \tilde V \triangleq V - \hat V \vspace{-0.cm}
\]
then the observer's objective is the exponential stability of $(\tilde R,\tilde V)=(I_3,0)$ or $(\bar R,\tilde V)=(I_3,0)$. These two possibilities are studied next.

From \eqref{attitudeSystem}, \eqref{dotv} and \eqref{observerGen}, one verifies that the dynamics of $(\tilde R, \tilde V)$ satisfy \vspace{-0.cm}
\begin{equation}\label{observerErrorDyn1}
\left\{
\begin{array}{ll}
&\!\!\!\!\!\!\!\! \dot{\tilde R} = \tilde R [\hat R \sigma_R]_\times \\
&\!\!\!\!\!\!\!\! \dot{\tilde V} = -[\Omega]_\times \tilde V + g \hat R^\top (\tilde R^\top - I_3) e_3 + \sigma_V
\end{array}\right. \vspace{-0.cm}
\end{equation}
while the dynamics of the error $(\bar R, \tilde V)$ are \vspace{-0.cm}
\begin{equation}\label{observerErrorDyn2}
\left\{
\begin{array}{ll}
&\!\!\!\!\!\!\!\! \dot{\bar R} = \bar R [\Omega]_\times  - [\Omega - \sigma_R]_\times \bar R\\
&\!\!\!\!\!\!\!\! \dot{\tilde V} = -[\Omega]_\times \tilde V + g (\bar R^\top - I_3)\hat R^\top e_3 + \sigma_V
\end{array}\right. \vspace{-0.cm}
\end{equation}
The next step consists in working out first order approximations of the error systems \eqref{observerErrorDyn1} and \eqref{observerErrorDyn2} complemented with first order approximations of the measurement equations. The application to these approximations of the Riccati observer design framework reported in \cite{HamSamTAC16} (see Section \ref{subsec:Riccati}) will then provide us with the equations of the proposed observers.

For this application we need to make the following technical (but non-restrictive) assumption.
\vspace{-0.2cm}
\begin{hypo}\label{assumpVOmbound}
The vehicle's velocities $v(t)$ and $\Omega(t)$ are bounded in norm by some positive numbers $v_{max}$ and $\Omega_{max}$, i.e. $|v(t)| \leq v_{max}$ and $|\Omega(t)| \leq \Omega_{max}$.
\end{hypo}

\vspace{-0.cm}
First order approximations of the attitude error equations are derived by considering (local) minimal parametrizations of the three-dimensional group of rotations $\mathrm{SO(3)}$. The parametrizations here chosen are the vector part $\tilde q$ (resp. $\bar q$) of the Rodrigues unit quaternion $\tilde Q = (\tilde q_0, \tilde q)$ (resp. $\bar Q = (\bar q_0, \bar q)$) associated with $\tilde R$ (resp. $\bar R$).
Rodrigues formula relating $\tilde Q$ (resp. $\bar Q$) to $\tilde R$ (reps. $\bar R$) are \vspace{-0.1cm}
\[
\begin{array}{l}
\tilde R = I_3 + 2 [\tilde q]_\times ( \tilde q_0 I_3 + [\tilde q]_\times)  \\
\bar R = I_3 + 2 [\bar q]_\times ( \bar q_0 I_3 + [\bar q]_\times)
\end{array} \vspace{-0.1cm}
\]
From these relations, and using also the fact that $\tilde q_0 = \pm\sqrt{1 - |\tilde q|^2}$ by definition of a unit quaternion, one deduces  \vspace{-0.1cm}
\[
\begin{array}{l}
\tilde R = I_3 + [\tilde \lambda]_\times + O(|\tilde \lambda|^2) ,\quad \mathrm{with} \,\, \tilde \lambda \triangleq 2 \tilde q  \\
\bar R = I_3 + [\bar \lambda]_\times  + O(|\bar \lambda|^2),\quad \mathrm{with} \,\, \bar \lambda \triangleq 2 \bar q
\end{array} \vspace{-0.1cm}
\]
Then, in view of the dynamics of $\tilde R$ and $\bar R$ in \eqref{observerErrorDyn1} and \eqref{observerErrorDyn2} one verifies (see also \cite{HamSamTAC16}) that the time-variations of $\tilde \lambda$ and $\bar \lambda$ satisfy the following equations \vspace{-0.cm}
\[
\begin{array}{l}
\dot{\tilde \lambda} = \hat R \sigma_R + O(|\tilde \lambda| |\sigma_R|)\\
\dot{\bar \lambda} = -[\Omega]_\times \bar \lambda + \sigma_R+ O(|\Omega||\bar\lambda|^2) + O(|\bar\lambda| |\sigma_R|)
\end{array}\vspace{-0.cm}
\]
As for the dynamics of $\tilde V$ one obtains, depending on the parametrization $\lambda$ or $\bar \lambda$ used for the attitude error \vspace{-0.cm}
\[
\begin{array}{l}
\dot{\tilde V} = -[\Omega]_\times \tilde V + g \hat R^\top[e_3]_\times \tilde \lambda  + \sigma_V +  O(|\tilde \lambda|^2)\\
\dot{\tilde V} = -[\Omega]_\times \tilde V  + g [\hat R^\top e_3]_\times \bar \lambda + \sigma_V +  O(|\bar \lambda|^2)
\end{array} \vspace{-0.cm}
\]

Concerning the measurement of $v_3$, in combination with the use of $\tilde \lambda$, one has \vspace{-0.cm}
\[
\begin{split}
&v_3 - e_3^\top \hat R \hat V  = e_3^\top R V - e_3^\top \hat R \hat V \\
&\quad = e_3^\top(\tilde R - I_3) \hat R (\hat V + \tilde V) + e_3^\top \hat R\tilde V\\
&\quad = -e_3^\top [\hat R \hat V]_\times \tilde \lambda + e_3^\top \hat R \tilde V + O(|\tilde \lambda||\tilde V|) + O(|V||\tilde \lambda|^2)
\end{split}
\]
and, in combination with the use of $\bar \lambda$ \vspace{-0.cm}
\[
\begin{split}
&v_3 - e_3^\top \hat R \hat V 
 = e_3^\top\hat R (\bar R - I_3) (\hat V + \tilde V) + e_3^\top \hat R\tilde V\\
&\quad = -e_3^\top \hat R [\hat V]_\times \bar \lambda + e_3^\top \hat R \tilde V
+ O(|\bar \lambda||\tilde V|) + O(|V||\bar \lambda|^2)
\end{split}
\]
As for the measurement of $m_{\mathcal B}$ one obtains respectively \vspace{-0.cm}
\[
\begin{array}{l}
\hat R m_{\mathcal B} - m_{\mathcal I} = (\tilde R^\top - I_3) m_{\mathcal I} = [m_{\mathcal I}]_\times \tilde \lambda + O(|\tilde\lambda|^2)\\
\hat R m_{\mathcal B} - m_{\mathcal I} =
\hat R (\bar R^\top - I_3) \hat R^\top m_{\mathcal I} = [m_{\mathcal I}]_\times \hat R\bar \lambda + O(|\bar\lambda|^2)
\end{array}\vspace{-0.cm}
\]
Note that one may also use the approximation $\hat R m_{\mathcal B} \times m_{\mathcal I} \approx \pi_{m_{\mathcal I}} \tilde \lambda$ when using the parametrization $\tilde \lambda$ for the attitude error.

\vspace{0.2cm}
In view of the previous relations, by setting the system output vector equal to \vspace{-0.1cm}
\begin{equation}\label{Y}
y = \begin{bmatrix} V_1 - \hat V_1 \\ V_2 - \hat V_2 \\ v_3 - e_3^\top \hat R \hat V \\
\hat R m_{\mathcal B} - m_{\mathcal I}
\end{bmatrix} \vspace{-0.1cm}
\end{equation}
one obtains LTV first order approximations in the form \eqref{LTV} with \vspace{-0.1cm}
\begin{equation}\label{paramtildelambda}
\!\left\{
\begin{split}
& x =\begin{bmatrix} \tilde \lambda \\ \tilde V \end{bmatrix}, \,\,x_1 = \tilde \lambda, \,\, x_2 = \tilde V, \,\,
u = \begin{bmatrix}\hat R \sigma_R \\ \sigma_V \end{bmatrix}, \\
&A \!=\! \begin{bmatrix} 0_{3\times 3} & \!\!0_{3\times 3}\\ g \hat R^\top[e_3]_\times & \!\!-[\Omega]_\times \end{bmatrix}\!, \, C \!=\! \begin{bmatrix} 0_{1\times 3} & \!\!e_1^\top  \\ 0_{1\times 3}  & \!\!e_2^\top \\
-e_3^\top [\hat R \hat V]_\times & \!\!e_3^\top \hat R \\
 [m_{\mathcal I}]_\times & \!\!0_{3\times 3} \end{bmatrix}
\end{split}\right. \!\!\!\!\!
\end{equation}
when using the parametrization $\tilde \lambda$, and \vspace{-0.1cm}
\begin{equation}\label{parambarlambda}
\!\left\{
\begin{split}
& x =\begin{bmatrix} \bar \lambda \\ \tilde V \end{bmatrix}, \,\,x_1 = \bar \lambda, \,\, x_2 = \tilde V, \,\,
u = \begin{bmatrix}\sigma_R \\ \sigma_V \end{bmatrix}, \\
&A \!=\! \begin{bmatrix} -[\Omega]_\times & \!\!0_{3\times 3}\\ g [\hat R^\top e_3]_\times & \!\!-[\Omega]_\times \end{bmatrix}\!, \, C \!=\! \begin{bmatrix} 0_{1\times 3} & \!\!e_1^\top  \\ 0_{1\times 3}  & \!\!e_2^\top \\
-e_3^\top \hat R [\hat V]_\times & \!\!e_3^\top \hat R \\
[m_{\mathcal I}]_\times \hat R &\!\! 0_{3\times 3} \end{bmatrix}
\end{split}\right. \!\!\!\!\!
\end{equation}
when using the parametrization $\bar \lambda$.

From there the observer associated with either one of the attitude error parametrizations is given by \eqref{observerGen} with $\sigma_R$ and $\sigma_V$ determined from the input $u$ calculated according to \eqref{U} and \eqref{CRE}.

\section{Observability analysis} \label{sec:observability_analysis}
According to \cite[Corollary 3.2]{HamSamTAC16}, good conditioning of the solutions $P(t)$ to the CREs and exponential stability of the previously derived observers rely on the uniform observability of the pair $(A^\star(t), C^\star(t))$ obtained by setting $x=0$ in the expressions of the matrices $A$ and $C$ derived previously.

In view of \eqref{paramtildelambda} one has  
\[
A^\star \!=\! \begin{bmatrix} 0_{3\times 3} & \!\!0_{3\times 3}\\[1ex] g R^\top[e_3]_\times & \!\!-[\Omega]_\times \end{bmatrix}, \,\,
C^\star \!=\! \begin{bmatrix} \begin{bmatrix} 0_{2\times 3} \\[1ex] -e_3^\top [v]_\times \end{bmatrix} & \Delta \\
 [m_{\mathcal I}]_\times & \!\!0_{3\times 3} \end{bmatrix} 
\]
when using $\tilde \lambda$ and, in view of \eqref{parambarlambda} 
\[
A^\star \!=\! \begin{bmatrix} -[\Omega]_\times & \!\!0_{3\times 3}\\[1ex] g[R^\top e_3]_\times & \!\!-[\Omega]_\times \end{bmatrix}, \,\,
C^\star \!=\! \begin{bmatrix} \begin{bmatrix} 0_{2\times 3} \\[1ex] -e_3^\top [v]_\times R \end{bmatrix} & \Delta \\
 [m_{\mathcal I}]_\times R & \!\!0_{3\times 3} \end{bmatrix} 
\]
when using $\bar \lambda$.

Define $D \!=\! \begin{bmatrix} D_{1,1} & D_{1,2} \\ D_{2,1} & D_{2,2}\end{bmatrix} \!\triangleq\! M^\top M$, with $M \triangleq \begin{bmatrix} C^\star \\ C^\star A^\star + \dot C^\star\end{bmatrix}$.
From Lemma \ref{lemObservabilityUniform} the pair $(A^\star, C^\star)$ is uniformly observable if $\exists \delta >0$, $\mu >0$ such that 
\begin{equation}\label{UOcondition}
\frac{1}{\delta} \int_{t}^{t+\delta} D(s) ds \geq \mu I_6 ,\,\, \forall t>0 
\end{equation}
Straightforward calculations yield 
\begin{equation}\label{Dlambdatilde}
\left\{
\begin{split}
D_{1,1}& =  \pi_{m_{\mathcal I}} + g^2 \pi_{e_3} - g^2(e_3 \times R e_3)(e_3 \times R e_3)^\top \\
&\quad + (e_3\times v)(e_3\times v)^\top + (e_3\times \dot v)(e_3\times \dot v)^\top \\
D_{1,2}& = D_{2,1}^\top = [v]_\times e_3 e_3^\top R + g[e_3]_\times R \pi_{e_3}[\Omega]_\times \\
D_{2,2}& = \Delta^\top \Delta -[\Omega]_\times\pi_{e_3}[\Omega]_\times
\end{split}\right.
\end{equation}
in the case of the $\tilde \lambda$ parametrization, and
\begin{equation}\label{Dlambdabar}
\left\{
\begin{split}
D_{1,1}& =  R^\top \big( \pi_{m_{\mathcal I}} + g^2 \pi_{e_3} - g^2(e_3 \times R e_3)(e_3 \times R e_3)^\top \\
&\quad + (e_3\times v)(e_3\times v)^\top + (e_3\times \dot v)(e_3\times \dot v)^\top \big)R \\
D_{1,2}& = D_{2,1}^\top = R^\top ([v]_\times e_3 e_3^\top R + g[ e_3]_\times R \pi_{e_3}[\Omega]_\times ) \\
D_{2,2}& = \Delta^\top \Delta -[\Omega]_\times\pi_{e_3}[\Omega]_\times
\end{split}\right.
\end{equation}
in the case of the $\bar \lambda$ parametrization. The determination of more explicit conditions whose satisfaction ensures uniform observability (and thus local exponential stability of the proposed observers) is not an easy task. However, it is possible to work out particular cases for which the stronger condition
\begin{equation}\label{IOUOCondition}
D(t) \geq \mu I_6,\,\, \forall t\geq 0 
\end{equation}
is satisfied. Both instantaneous observability and uniform observability of the pair $(A^\star, C^\star)$ are then granted (see Remark \ref{rem1}). The following lemma points out such particular cases. Its proof is based on the fact that \eqref{IOUOCondition} is equivalent to
\begin{equation*}\label{IOUOConditionBIS}
X^\top D(t) X \geq \mu |X|^2, \quad\forall X  \in \mathbb{R}^6, \forall t\geq 0  
\end{equation*}
which, by simple computations, is also equivalent to $\forall x,y \in \mathbb{R}^3$ and $\forall t\geq 0$ 
\begin{equation}\label{IOUOConditionTER}
\begin{split}
&|m_{\mathcal{I}} \!\times\! x|^2\! +\! g^2\big|[Re_3]_\times [e_3]_\times x\big|^2 \\
&\quad\,\,\, + \left((e_3 \!\times\! v)^\top \! x \right)^2 \!+\! \left((e_3 \!\times\! \dot v)^\top \!x \right)^2 \\
&\quad\,\,\, +y_1^2 + y_2^2 + \left((R^\top e_3)^\top y\right)^2 + \big|[e_3]_\times [\Omega]_\times y \big|^2 \\
&\quad\,\,\, + 2g(R^\top [Re_3]_\times [e_3]_\times x)^\top ([e_3]_\times [\Omega]_\times y)\\
&\quad\,\,\, + 2 ((v \!\times\!  e_3)^\top \!x) \left((R^\top e_3)^\top y\right) \;\geq\; \mu (|x|^2 + |y|^2)
\end{split}
\end{equation}

\begin{lem}\label{IOLemma} Assume that \vspace{-.cm}
\begin{equation}\label{IOConditionCases}
\exists \rho > 0 \,\,\,\, \mathrm{s.t.} \,\,\,\,  |R_{3,3}(t)| \geq \rho, \,\,\forall t\geq 0 \vspace{-.cm}
\end{equation}
Then, condition \eqref{IOUOCondition} is satisfied in the following cases:
\begin{enumerate}
\item Motion along the vertical direction, i.e. $v(t) \times e_3  \equiv 0$;
\item Pure translation, i.e. $\Omega(t) \equiv 0$;
\item Slow motion such that $v_{max} \Omega_{max} \leq \frac{g\rho^2}{\sqrt{6}}$, with $v_{max}$ and $\Omega_{max}$ standing for the bounds of $v$ and $\Omega$ defined in Assumption \ref{assumpVOmbound}.
\end{enumerate}
\end{lem}

\vspace{-.cm}
(Proof in Appendix \ref{appendix1}).

\vspace{.2cm}
Some comments about condition \eqref{IOConditionCases} are in order. This condition indicates that the gravity direction expressed in the body-fixed frame $\{{\mathcal B}\}$, i.e. $R^{\top} e_3$, never crosses the plane spanned by $e_1$ and $e_2$ or approaches it asymptotically.
For instance, if $\forall t:R(t)^\top e_3 \in span{(e_1,e_2)}$ then both observability conditions  \eqref{UOcondition} and \eqref{IOUOCondition} are not satisfied since in that case the last row and last column of $D$ (given by \eqref{Dlambdatilde}) are equal to zero. However, this very particular situation of non observability is not supposed to occur in the case of quadrotor UAV navigation. Now $R^{\top} e_3$ may temporarily cross the plane $span{(e_1,e_2)}$ thus leading to the violation of condition \eqref{IOConditionCases} as well as of the instantaneous observability condition \eqref{IOUOCondition}. But the condition \eqref{UOcondition} of uniform observability may still be satisfied in this case.

\vspace{-.1cm}
\begin{lem}\label{UOLemma} Assume that \vspace{-.cm}
\begin{equation}\label{UOConditionCases}
\exists \delta, \rho > 0 \,\,\,\, \mathrm{s.t.} \,\,\,\,  \frac{1}{\delta} \int_t^{t+\delta}   R_{3,3}^2(s)ds \geq \rho^2, \,\,\forall t\geq 0 \vspace{-.cm}
\end{equation}
Then, condition \eqref{UOcondition} is satisfied in the same three cases as in Lemma \ref{IOLemma} and also in the case of
persistently accelerated translational motion such that $\exists \bar\delta, \bar\rho > 0$, $\forall t\geq 0 $:
\begin{equation}\label{UOConditionDotV}
\min\left(\frac{1}{\bar\delta} \int_t^{t+\bar\delta}   \dot v_1^2(s)ds, \frac{1}{\bar\delta} \int_t^{t+\bar\delta}   \dot v_2^2(s)ds \right)\geq \bar\rho \vspace{-.cm}
\end{equation}
\end{lem}

(Proof in Appendix \ref{appendix2}).

\vspace{-.cm}
\begin{rem}\label{remarkMagAbsence}
In the case where magnetometer measurements are absent the observability condition \eqref{UOcondition} is never satisfied.
Indeed, by inspection of the expression of $D_{1,1}$ (given in relation \eqref{Dlambdatilde}) from which the term $\pi_{m_{\mathcal I}}$, is removed
one easily verifies that the third row and third column of this matrix are equal to zero. One also verifies that the third row and third column of $D_{1,2}$ (and $D_{2,1}$) are equal to zero. Therefore, the third row and third column of $D$ are equal to zero. This clearly forbids the satisfaction of the condition \eqref{UOcondition}.
\end{rem}

\section{Simulation results}\label{sec:simulation}

\begin{figure}[!t]\centering%
\vspace*{-0.cm}
\psfrag{Real value}{\tiny Real value}%
\psfrag{Observer 1}{\tiny Observer 1}%
\psfrag{Observer 2}{\tiny Observer 2}%
\psfrag{t (s)}{\scriptsize $t (s)$}%
\psfrag{roll (deg)}{\scriptsize $\text{roll}\, \phi,\hat\phi$}%
\psfrag{pitch (deg)}{\scriptsize $\text{pitch}\, \theta,\hat\theta$}%
\psfrag{yaw (deg)}{\scriptsize $\text{yaw}\, \psi,\hat\psi$}%
\includegraphics[width=1.\linewidth, height = 8cm]{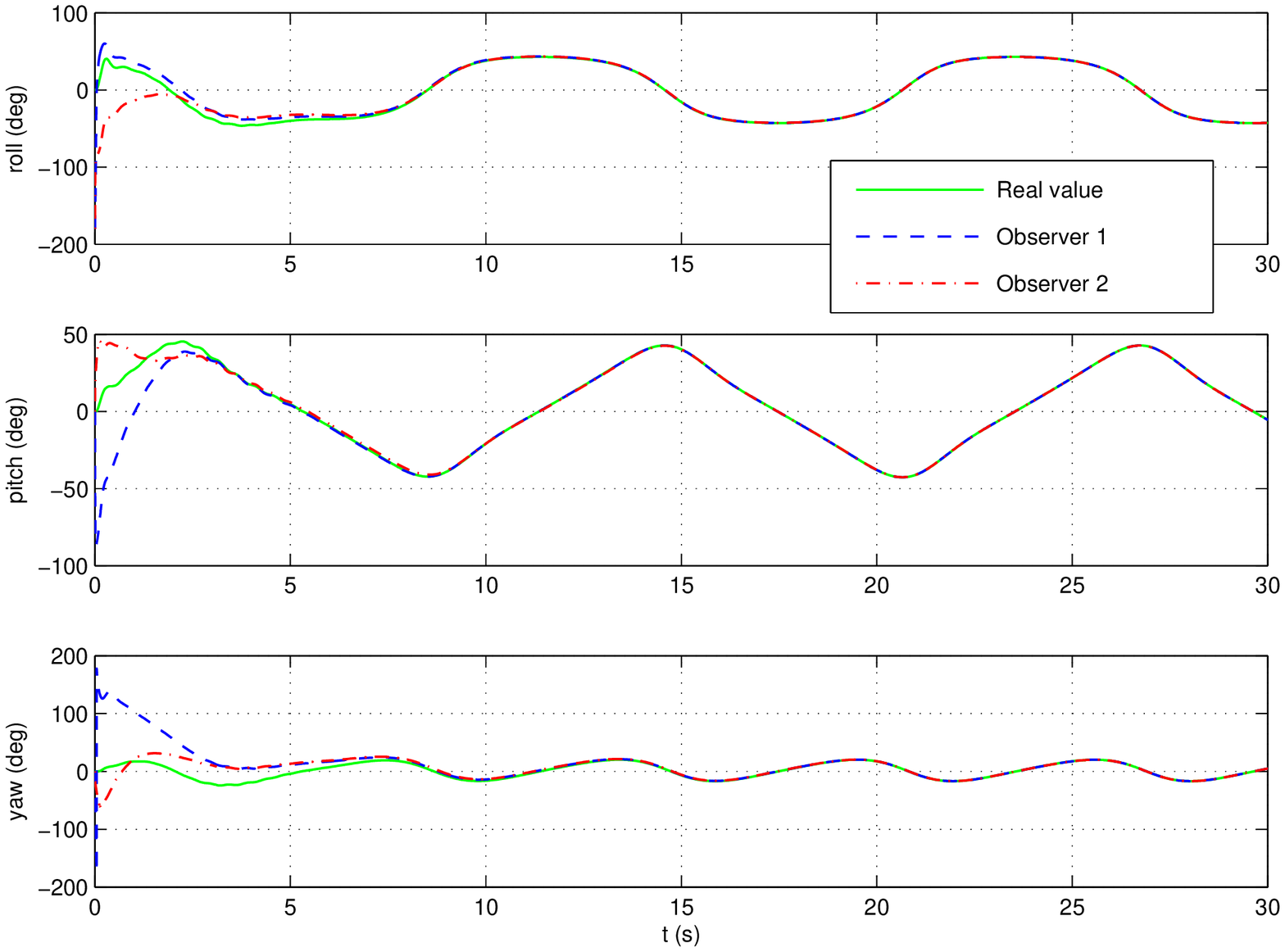}\vspace*{-0.6cm}
\caption{Simulation 1: Estimated and real attitude represented by roll, pitch and yaw Euler angles ($deg$) versus time ($s$).} \label{figsimu1a}
\vspace*{-0.cm}
\psfrag{Real value}{\tiny Real value}%
\psfrag{Observer 1}{\tiny Observer 1}%
\psfrag{Observer 2}{\tiny Observer 2}%
\psfrag{t (s)}{\scriptsize $t (s)$}%
\psfrag{v1 (m/s)}{\scriptsize $V_1$, $\hat V_1$}%
\psfrag{v2 (m/s)}{\scriptsize $V_2$, $\hat V_2$}%
\psfrag{v3 (m/s)}{\scriptsize $V_3$, $\hat V_3$}%
\includegraphics[width=1.\linewidth, height = 8cm]{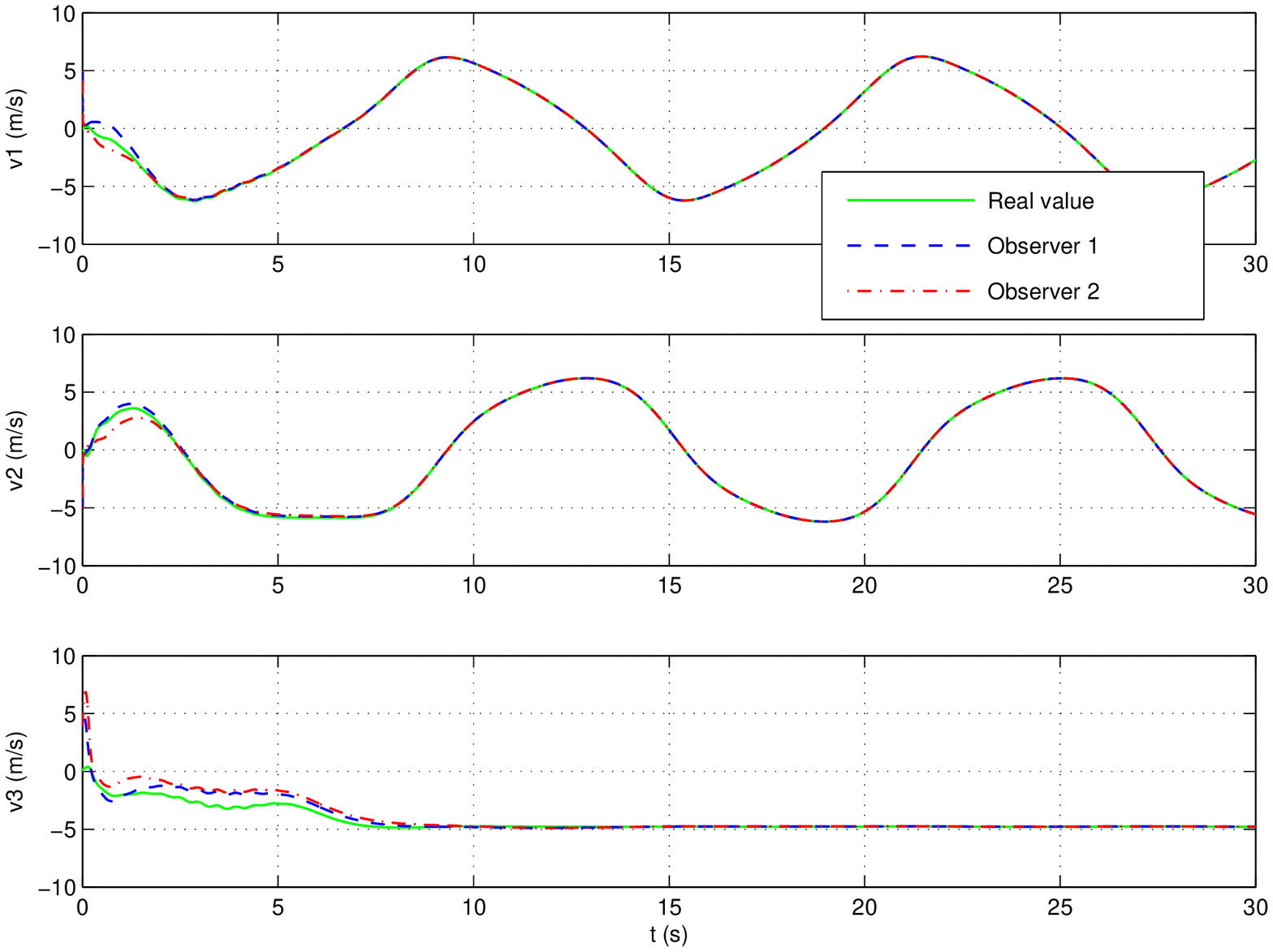} \vspace*{-1.cm}
\caption{Simulation 1: Estimated and real velocity ($m/s$) versus time ($s$), where $V_i$ (resp. $\hat V_i$) is the $i^{th}$ component of $V$ (resp. $\hat V$).} \label{figsimu1b}  \vspace*{-0.2cm}
\end{figure}

\begin{figure}[!ht]\centering%
\vspace*{-0.cm}
\psfrag{Real value}{\tiny Real value}%
\psfrag{Observer 1}{\tiny Observer 1}%
\psfrag{Observer 2}{\tiny Observer 2}%
\psfrag{t (s)}{\scriptsize $t (s)$}%
\psfrag{roll (deg)}{\scriptsize $\text{roll}\, \phi,\hat\phi$}%
\psfrag{pitch (deg)}{\scriptsize $\text{pitch}\, \theta,\hat\theta$}%
\psfrag{yaw (deg)}{\scriptsize $\text{yaw}\, \psi,\hat\psi$}%
\includegraphics[width=1.\linewidth, height = 8cm]{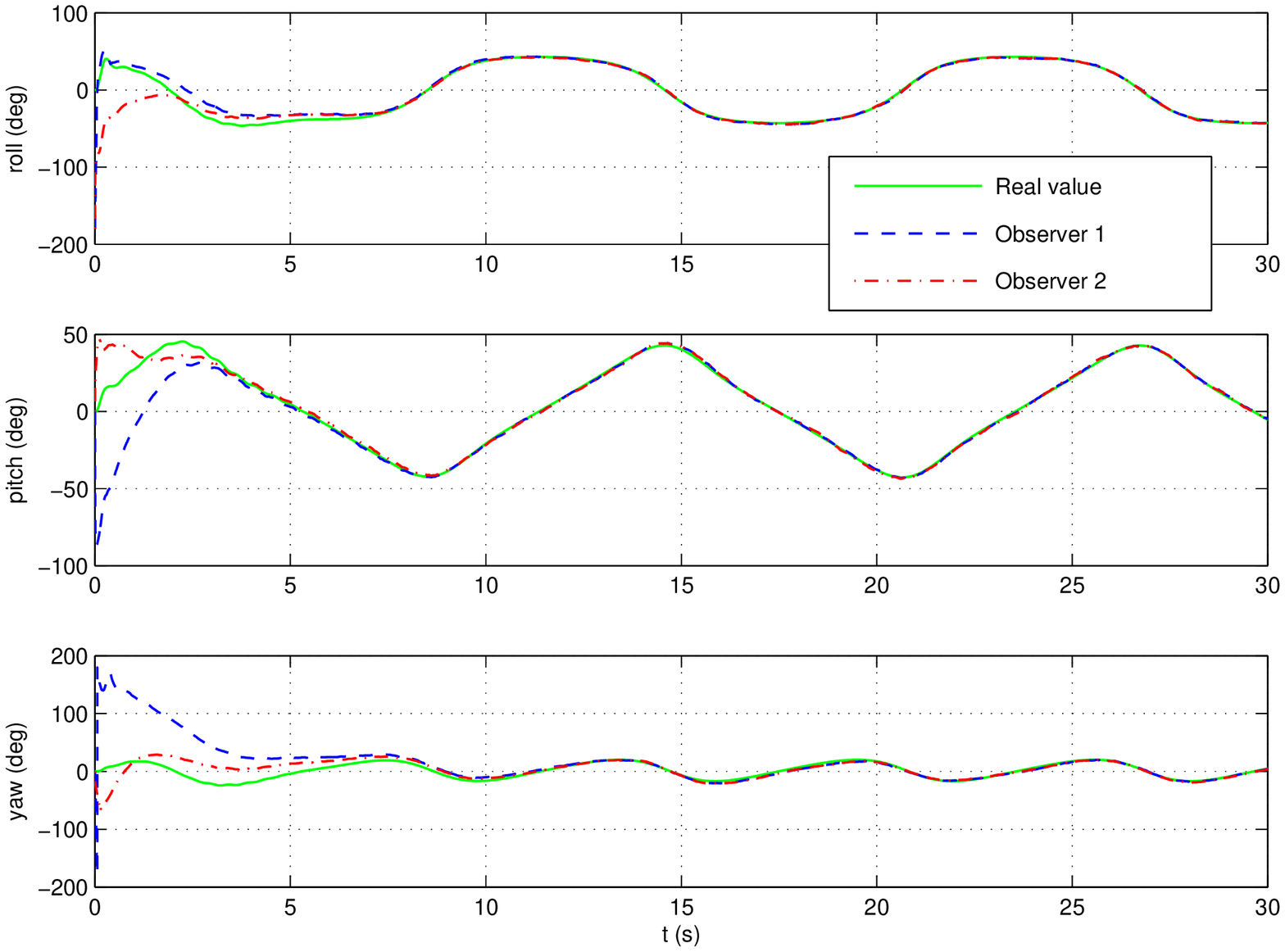}\vspace*{-0.6cm}
\caption{Simulation 2: Estimated and real attitude represented by roll, pitch and yaw Euler angles ($deg$) versus time ($s$).} \label{figsimu2a}
\psfrag{Real value}{\tiny Real value}%
\psfrag{Observer 1}{\tiny Observer 1}%
\psfrag{Observer 2}{\tiny Observer 2}%
\psfrag{t (s)}{\scriptsize $t (s)$}%
\psfrag{v1 (m/s)}{\scriptsize $V_1$, $\hat V_1$}%
\psfrag{v2 (m/s)}{\scriptsize $V_2$, $\hat V_2$}%
\psfrag{v3 (m/s)}{\scriptsize $V_3$, $\hat V_3$}%
\includegraphics[width=1.\linewidth, height = 8cm]{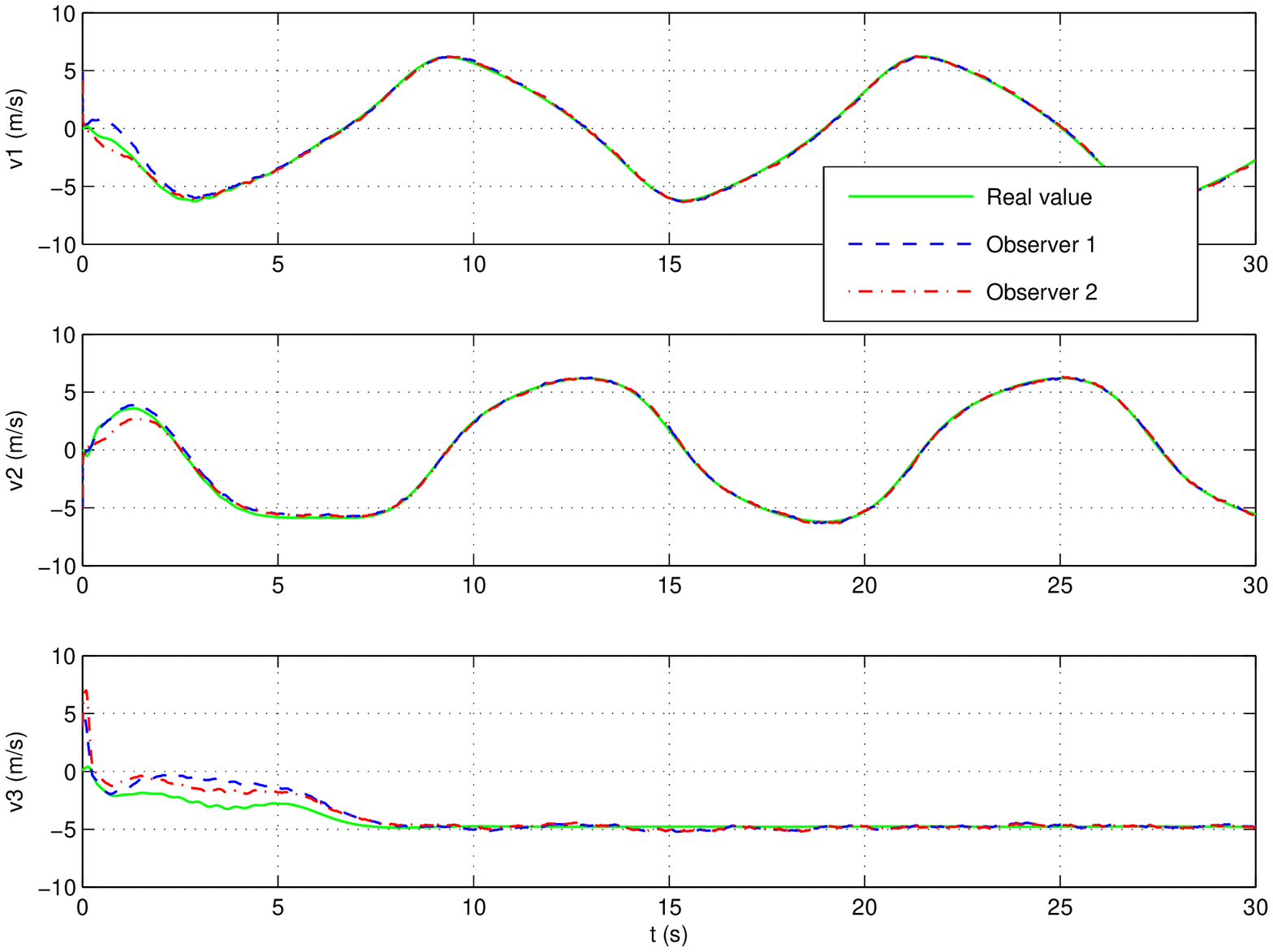} \vspace*{-1.cm}
\caption{Simulation 2: Estimated and real velocity ($m/s$) versus time ($s$).} \label{figsimu2b}  \vspace*{-0.cm}
\end{figure}

Simulations are conducted on a model of a vertical take-off and landing (VTOL) aerial drone, also used in \cite{hua10cep}. The vehicle is stabilised along a circular reference trajectory, with the linear velocity expressed in the inertial frame $\{\mathcal I\}$ given by $v_r =[-15\alpha\sin(\alpha t);15\alpha\cos(\alpha t);0]$ $(m/s)$, with $\alpha=2/\sqrt{15}$. Due to aerodynamic forces the vehicle's orientation varies in large proportions.
The normalized earth's magnetic field and the gravity constant are respectively equal to $m_{\mathcal I}= [0.434; -0.0091; 0.9008]$ and $g=9.81 (m/s^2)$.

In the absence of other works addressing the same problem of coupled velocity-aided attitude estimation, comparisons are only carried out for the two observers here proposed. We call  {\em Observer 1} (resp. {\em Observer 2}) the observer derived with the parametrization $\tilde \lambda$ (resp. $\bar \lambda$).

Both observers are tuned analogously to Kalman-Bucy filters where the matrices $S$ and $Q^{-1}$ are interpreted as covariance matrices of the additive noise on the system state and output respectively. The following parameters
are chosen for both observers: \vspace{-0.1cm}
\[
\left\{
\begin{array}{l}
P(0)=\mathrm{diag}(2 I_3, 20I_3)\\
Q(t)=\mathrm{diag}(25I_3, 100I_3)\\
S(t)=\mathrm{diag}(0.01I_3, I_3)
\end{array}\right. \vspace{-0.1cm}
\]

Two simulations are reported hereafter.

\vspace{0.1cm}
\noindent $\bullet$ {\bf Simulation 1:} In this simulation, the observers are simulated in the ideal case (i.e. noise-free measurements) for a set of initial attitude and velocity estimates corresponding to the following initial estimation errors
\begin{equation}\label{simu1init}
\left\{\!\!\!\!\!\!\!\!
\begin{array}{ll}
&\tilde v(0)= [-5;5;-5] (m/s) \\[1ex]
&\tilde q(0) = [\cos(\frac{\pi}{2}); \sin(\frac{\pi}{2}) e_1]
\end{array}\right.
\end{equation}
This extreme case corresponds to an initial attitude error of $180(deg)$ in roll w.r.t. the true attitude. The time evolutions of the estimated and real attitudes, represented by Euler angles, along with the estimated and real velocity are shown in Figs. \ref{figsimu1a} and  \ref{figsimu1b}, respectively. Both observers ensure the asymptotic convergence of the estimated variables to the real values despite the extremely large initial estimation errors. Their convergence rates are similar and quite satisfactory.

\vspace{0.2cm}
\noindent $\bullet$ {\bf Simulation 2:} This simulation is conducted with the same initial condition \eqref{simu1init} as in Simulation 1. However, the measurements are now corrupted by Gaussian zero-mean additive noises with standard deviations reflecting the above choice of $Q$ ($0.2\,m/s$ for $v_3$ and $V_{1,2}$ and $0.1$ for $m_{\mathcal B}$) and of $S$ ($0.1 \,rad/s$ for $\Omega$ and $1\,m/s^2$ for $a_{\cal B}$). Moreover, they are discretized with update frequencies of $20\, Hz$, for the measurements of $V_{1}$, $V_2$, $v_3$ and $m_{\mathcal B}$, and of $50 \, Hz$, for the measurements of $\Omega$ and $a_{\mathcal B}$. The results reported in Figs. \ref{figsimu2a} and \ref{figsimu2b} indicate that important noises and low update frequencies of the measurements only marginally affect the overall performance of both observers.

\section{Conclusions}\label{sec:conclusions}
In this paper, a new problem of coupled velocity-aided attitude estimation has been addressed and two nonlinear observers have been proposed on the basis of a recent deterministic Riccati observer design framework. They are supported by comprehensive stability and observability analysis, and also by convincing simulation results.

\section*{Acknowledgement}
This research was supported by the French {\em Agence Nationale de la Recherche} via the ASTRID SCAR (ANR-12-ASTR-0033) and ROBOTEX (ANR-10-EQPX-44) projects.

\appendix
\setcounter{section}{0}
\subsection{Proof of Lemma \ref{IOLemma}} \label{appendix1}
\noindent {\bf Case 1}:
In this case $v \times e_3\equiv \dot v \times e_3 \equiv 0$ and \vspace{-.1cm}
\begin{equation}\label{IOUOConditionTER1}
\begin{split}
&\!\!\!\!X^\top D X = |m_{\mathcal{I}} \!\times\! x|^2 \! +\!\frac{\varepsilon_r g^2}{1+\varepsilon_r}\big|[Re_3]_\times [e_3]_\times x\big|^2   \\
& \quad\!\!\!\!+y_1^2 + y_2^2 + \left((R^\top e_3)^\top y\right)^2 -  \varepsilon_r \big|[e_3]_\times [\Omega]_\times y \big|^2 \\
&\quad \!\!\!\! +\!\left(\frac{g}{\sqrt{1\!+\!\varepsilon_r}} R^\top [Re_3]_\times [e_3]_\times x  + \sqrt{1\!+\!\varepsilon_r} [e_3]_\times [\Omega]_\times y\right)^2  \\
& \!\!\!\! \geq\, |m_{\mathcal{I}} \!\times\! x|^2 \! +\!\frac{\varepsilon_r g^2}{1+\varepsilon_r}\big|[Re_3]_\times [e_3]_\times x\big|^2   \\
&\quad\!\!\!\!  +y_1^2 + y_2^2 + \left((R^\top e_3)^\top y\right)^2 -  \varepsilon_r \Omega_{max}^2 |y|^2
\end{split}
\end{equation}

\vspace{-0.2cm}
\noindent with $\varepsilon_r>0$  such that $\exists \mu_y^{r} >0$:  \vspace{-.1cm}
\begin{equation}\label{findepsPureRotation}
y_1^2 + y_2^2 + \left((R^\top e_3)^\top y\right)^2 -  \varepsilon_r \Omega_{max}^2 |y|^2  \geq \mu_y^{r} |y|^2  \vspace{-.1cm}
\end{equation}
A number $\varepsilon_r$ satisfying this inequality is calculated next. Defining $\gamma \triangleq R^\top e_3 \in S^2$ one gets \vspace{-.1cm}
\[
\begin{array}{ll}
&\!\!\!\!\!\!\left(\gamma^\top y\right)^2 = (\gamma_3 y_3)^2 \!+\! 2 (\gamma_1 y_1 \!+\! \gamma_2 y_2)(\gamma_3 y_3) \!+ \!(\gamma_1 y_1 \!+ \!\gamma_2 y_2)^2 \\[1ex]
&\geq \frac{1}{3} (\gamma_3 y_3)^2 - \frac{1}{2}(\gamma_1 y_1 + \gamma_2 y_2)^2 \\[1ex]
&\geq \frac{1}{3} (\gamma_3 y_3)^2 - (\gamma_1 y_1)^2 - ( \gamma_2 y_2)^2  \vspace{-.cm}
\end{array}
\]
when using the following Young inequalities   \vspace{-.cm}
\[
\begin{array}{ll}
& 2 (\gamma_1 y_1 + \gamma_2 y_2)(\gamma_3 y_3) \geq -\frac{2}{3} (\gamma_3 y_3)^2  -\frac{3}{2} (\gamma_1 y_1 + \gamma_2 y_2)^2\\[1ex]
& (\gamma_1 y_1 + \gamma_2 y_2)^2 \leq 2((\gamma_1 y_1)^2 + ( \gamma_2 y_2)^2)
\end{array}  \vspace{-.cm}
\]
Since $\gamma_3 = R_{3,3}$ one deduces from \eqref{IOConditionCases} that  \vspace{-.1cm}
\begin{equation}\label{ypartPureRota}
\begin{array}{ll}
&\!\!\!\!\!\!\!\!\!\!y_1^2 + y_2^2 + \left((R^\top e_3)^\top y\right)^2 -  \varepsilon_r \Omega_{max}^2 |y|^2  \\[1ex]
&\geq (1-\gamma_1^2) y_1^2 + (1-\gamma_2^2) y_2^2 + \frac{1}{3} \gamma_3^2 y_3^2 -  \varepsilon_r \Omega_{max}^2 |y|^2  \\[1ex]
&\geq (\frac{1}{3} \gamma_3^2 -  \varepsilon_r \Omega_{max}^2 )|y|^2 + \frac{2}{3} \gamma_3^2 (y_1^2+y_2^2)\\[1ex]
&\geq \mu_y^{r} |y|^2  \vspace{-.2cm}
\end{array}
\end{equation}
with $\mu_y^{r} \triangleq {\rho^2}/{3} -  \varepsilon_r \Omega_{max}^2$.
Therefore, any number $\varepsilon_r$ such that $0<\varepsilon_r< {\rho^2}/{(3\Omega_{max}^2})$ ensures  that $\mu_y^r$ in \eqref{ypartPureRota} is positive.

Let us now consider the term $|m_{\mathcal{I}} \!\times\! x|^2 \! +\!\frac{\varepsilon_r g^2}{1+\varepsilon_r}\big|[Re_3]_\times [e_3]_\times x\big|^2 $ involved in the last inequality of \eqref{IOUOConditionTER1}. By simple computations one obtains \vspace{-.1cm}
\[
[Re_3]_\times [e_3]_\times x \!=\! \begin{bmatrix} -R_{3,3} & 0& 0 \\ 0&-R_{3,3} &0 \\ R_{1,3} & R_{2,3} & 0\end{bmatrix} \!x \!=\!
\begin{bmatrix}  -R_{3,3} x_1 \\ -R_{3,3} x_2\\ R_{1,3} x_1 \!+ \!R_{2,3} x_2\end{bmatrix}\vspace{-.1cm}
\]
Thus, defining $\bar\varepsilon_r \triangleq \displaystyle\frac{\varepsilon_r g^2}{1+\varepsilon_r}$ and using \eqref{IOConditionCases} one deduces that\vspace{-.1cm}
\begin{equation}\label{xpartPureRota}
\begin{array}{ll}
&\!\!\!\!\!\!\!\! |m_{\mathcal{I}} \!\times\! x|^2 \! +\! \frac{\varepsilon_r g^2}{1+\varepsilon_r} \big|[Re_3]_\times [e_3]_\times x\big|^2 \\[1ex]
&\!\!\!\!= (m_{3} x_2 \!-\! m_{2} x_3)^2 \!+\! (m_{3} x_1 \!-\! m_{1} x_3)^2 \!+\! (m_{2} x_1 \!-\! m_{1} x_2)^2 \\
& + \bar\varepsilon_r R_{3,3}^2 (x_1^2+x_2^2) + \bar\varepsilon_r (R_{1,3} x_1 \!+ \!R_{2,3} x_2)^2 \\[1ex]
&\!\!\!\!\geq (m_{3} x_2 \!-\! m_{2} x_3)^2 \!+\! (m_{3} x_1 \!-\! m_{1} x_3)^2 + \bar\varepsilon_r \rho^2 (x_1^2+x_2^2)  \\[1ex]
&\!\!\!\!=  {\tiny\left(\! \sqrt{m_{3}^2\!+\!\frac{\bar\varepsilon_r\rho^2}{2}}x_2 \!-\! \frac{m_{2}m_3}{\sqrt{m_{3}^2 + \frac{\bar\varepsilon_r\rho^2}{2}}}  x_3 \!\right)^2} \\
& +  \left(\!\sqrt{m_{3}^2\!+\!\frac{\bar\varepsilon_r\rho^2}{2}} x_1 \!-\! \frac{m_1 m_3}{\sqrt{m_{3}^2 +\frac{\bar\varepsilon_r\rho^2}{2}}}  x_3 \!\right)^2 \\
& +\frac{\bar\varepsilon_r \rho^2}{2} (x_1^2+x_2^2)  + \frac{(m_1^2 + m_2^2) \bar\varepsilon_r\rho^2}{2m_{3}^2+ \bar\varepsilon_r\rho^2} x_3^2 \\[1ex]
&\!\!\!\!  \geq \mu_x^{r}  |x|^2
\end{array}
\vspace{-.1cm}
\end{equation}
with $\mu_x^{r} \triangleq \min\left(\frac{\bar\varepsilon_r \rho^2}{2}, \frac{(m_1^2 + m_2^2) \bar\varepsilon_r\rho^2}{2m_{3}^2+ \bar\varepsilon_r\rho^2}  \right)$. This number is positive since $m_{\mathcal{I}}$ and $e_3$ are non-collinear by assumption. From  \eqref{IOUOConditionTER1}, \eqref{ypartPureRota}, \eqref{xpartPureRota} one then deduces that $X^\top D X \geq \mu |X|^2$ with $\mu \triangleq \min(\mu_x^{r},\mu_y^{r})>0$. This concludes the proof of the first case.

\vspace{.2cm}
\noindent {\bf Case 2}:
The proof proceeds analogously to the proof of the first case. Since $\Omega(t) \equiv 0$, the left-hand side of \eqref{IOUOConditionTER} satisfies \vspace{-.1cm}
\begin{equation}\label{IOUOConditionTER2}
\begin{split}
\!\!\!\!X^\top D X &\geq y_1^2 + y_2^2 + \frac{\varepsilon_t}{1+\varepsilon_t}\left((R^\top e_3)^\top y\right)^2\\
&\!\!\!\!\!\!\!\!\!\!\!\!\!\!\!\! +\!|m_{\mathcal{I}} \!\times\! x|^2\! +\! g^2\big|[Re_3]_\times [e_3]_\times x\big|^2 \!-\!\varepsilon_t v_{max}^2 (x_1^2+ x_2^2)
\end{split}
\end{equation}
with $\varepsilon_t>0$ specified hereafter. Relation \eqref{xpartPureRota} is now replaced by \vspace{-.1cm}
\begin{equation}\label{xpartPureTrans}
\begin{split}
&\!\!\!\! |m_{\mathcal{I}} \!\times\! x|^2\! +\! g^2\big|[Re_3]_\times [e_3]_\times x\big|^2 -\varepsilon_t v_{max}^2 (x_1^2+x_2^2) \\
& \geq  \left(\frac{g^2 \rho^2}{2} -\varepsilon_t v_{max}^2\right) (x_1^2+x_2^2)  + \frac{(m_1^2 + m_2^2) g^2\rho^2}{2m_{3}^2+ g^2\rho^2} x_3^2 \\
& \geq \min\left(\frac{g^2 \rho^2}{2} -\varepsilon_t v_{max}^2, \frac{(m_1^2 + m_2^2) g^2\rho^2}{2m_{3}^2+ g^2\rho^2} \right) |x|^2 \\
& \geq \mu_x^{t} |x|^2
\end{split}\!\!\!\! \vspace{-.cm}
\end{equation}
with $\mu_x^{t} \triangleq \min\left(\frac{g^2 \rho^2}{2}- \varepsilon_t v_{max}^2, \frac{(m_1^2 + m_2^2) g^2\rho^2}{2m_{3}^2+ g^2\rho^2} \right)$. This number is positive if $\varepsilon_t$ is chosen such that
$0<\varepsilon_t <   \frac{g^2 \rho^2}{2 v_{max}^2}$.
Relation \eqref{ypartPureRota} is now replaced by  \vspace{-.1cm}
\begin{equation}\label{ypartPureTrans}
\begin{split}
  y_1^2 \!+\! y_2^2 &\!+\! \frac{\varepsilon_t}{1\!+\!\varepsilon_t}\left((R^\top\! e_3)^\top y\right)^2 \\
&\geq  \frac{\varepsilon_t}{1+\varepsilon_t} (y_1^2 + y_2^2 + \gamma^\top y) \geq \mu_y^{t}   |y|^2
\end{split}\!\!\!\!
\end{equation}

\vspace{-.2cm}
\noindent with $\mu_y^{t} \triangleq \displaystyle\frac{\varepsilon_t \rho^2}{3(1+\varepsilon_t)}$. From  \eqref{IOUOConditionTER2}, \eqref{xpartPureTrans}, \eqref{ypartPureTrans} one then deduces that $X^\top D(t) X \geq \mu |X|^2$ with $\mu \triangleq \min(\mu_x^{t},\mu_y^{t})>0$. This concludes the proof of the second case.

\vspace{.2cm}
\noindent {\bf Case 3}:
Using the same procedure as the one used to derive relations \eqref{IOUOConditionTER1}, \eqref{ypartPureRota}, and \eqref{xpartPureRota} one deduces that  \vspace{-.1cm}
\begin{equation*}\label{IOUOConditionTER3}
\begin{split}
&X^\top D X \geq \\
&|m_{\mathcal{I}} \!\times\! x|^2\! +\! \frac{\varepsilon_2 g^2}{1+\varepsilon_2}\big|[Re_3]_\times [e_3]_\times x\big|^2 \!-\!\varepsilon_1 v_{max}^2(x_1^2+x_2^2)\\
& +y_1^2 + y_2^2 + \frac{\varepsilon_1}{1+\varepsilon_1}\left((R^\top e_3)^\top y\right)^2 -  \varepsilon_2 \Omega_{max}^2 |y|^2 \\
&\geq \left(\frac{\varepsilon_2 g^2\rho^2}{2(1+\varepsilon_2)} \!-\!\varepsilon_1 v_{max}^2 \right)(x_1^2+x_2^2) +  \frac{(m_1^2 + m_2^2) \bar\varepsilon_2\rho^2}{2m_{3}^2+ \bar\varepsilon_2\rho^2} x_3^2 \\
&+\left( \frac{\varepsilon_1 \rho^2}{3(1+\varepsilon_1)} -  \varepsilon_2 \Omega_{max}^2\right) |y|^2
\end{split}
\end{equation*}

\vspace{-.3cm}
\noindent
with $\bar\varepsilon_2 \triangleq \displaystyle \frac{\varepsilon_2 g^2}{1+\varepsilon_2}$ and $\varepsilon_1, \varepsilon_2>0$ chosen such that \vspace{-.1cm}
\begin{equation}\label{eps12cond}
\left\{
\begin{array}{l}
\displaystyle\frac{\varepsilon_2}{1+\varepsilon_2} > \alpha_1 \varepsilon_1, \,\, \mathrm{with} \,\, \alpha_1 \triangleq \frac{2v_{max}^2}{g^2 \rho^2}\\[1ex]
\displaystyle\frac{\varepsilon_1}{1+\varepsilon_1} > \alpha_2 \varepsilon_2, \,\, \mathrm{with} \,\, \alpha_2 \triangleq \frac{3\Omega_{max}^2}{\rho^2}
\end{array}\right.
\end{equation}
One then verifies that positive solutions of $\varepsilon_1$ and $\varepsilon_2$ to \eqref{eps12cond} exist if $\alpha_1 \alpha_2 < 1$ or, equivalently, if
$v_{max} \Omega_{max} \leq \frac{g\rho^2}{\sqrt{6}}$. From here \eqref{IOUOCondition} follows immediately.

\subsection{Proof of Lemma \ref{UOLemma}} \label{appendix2}
Condition \eqref{UOcondition} is equivalent to the following: $\forall X=[x^\top,y^\top]^\top \in \mathbb{R}^6$, $\forall t\geq 0$ \vspace{-.1cm}
\[
\frac{1}{\delta} \int_t^{t+\delta}\!\!\!\!\!\!X^\top D(s)X ds \geq \mu |X|^2, \,\,\forall X=[x^\top\!,y^\top]^\top \!\!\in \mathbb{R}^6, \forall t\geq 0 \vspace{-.1cm}
\]
We only develop the proof for the first case. The proofs for the three other cases proceed analogously to the proof of the first case and the proof of Lemma \ref{IOLemma}.
Using the fact that
$v\times e_3 \equiv  \dot v\times e_3 \equiv 0$ one deduces that \vspace{-.1cm}
\begin{equation*}
\begin{array}{ll}
&\!\!\!\!\!\!\!\!\frac{1}{\delta} \int_t^{t+\delta}X^\top D(s)X ds   \\[1ex]
&\geq\frac{1}{\delta} \int_t^{t+\delta} (|m_{\mathcal{I}} \!\times\! x|^2 \! +\!\frac{\varepsilon_r g^2}{1+\varepsilon_r}\big|[R(s)e_3]_\times [e_3]_\times x\big|^2   \\[1ex]
&\quad +y_1^2 + y_2^2 + \left((R(s)^\top e_3)^\top y\right)^2 -  \varepsilon_r \Omega_{max}^2 |y|^2)ds  \\[1ex]
&\geq \min\left(\frac{\bar\varepsilon_r}{2}, \frac{(m_1^2 + m_2^2) \bar\varepsilon_r }{2m_{3}^2+ \bar\varepsilon_r} \right)
\left(\frac{1}{\delta} \int_t^{t+\delta}   R_{3,3}^2(s) ds \right) |x|^2  \\[1ex]
&+ \frac{1}{3} \left(\frac{1}{\delta} \int_t^{t+\delta}   R_{3,3}^2(s) ds - 3 \Omega_{max}^2\varepsilon_r\right) |y|^2
\end{array} \vspace{-.1cm}
\end{equation*}
with $\bar\varepsilon_r \triangleq \frac{\varepsilon_r g^2}{1+\varepsilon_r}$ and $0<\varepsilon_r < \frac{1}{3 \Omega_{max}^2 \delta} \int_t^{t+\delta}   R_{3,3}^2(s) ds$. Therefore, \eqref{UOcondition} holds if the ``persistent excitation'' condition  \eqref{UOConditionCases} is satisfied.


\begin{thebibliography}{99}
\bibitem{Allibert16} 
G. Allibert, R. Mahony, and M. Bangura. 
\newblock Velocity aided attitude estimation for aerial robotic vehicles using latent rotation scaling. 
\newblock In {\em IEEE Int. Conf. on Robotics and Automation}, pages 1538--1543, 2016.
%
\bibitem{batista2012}
P. Batista, C. Silvestre, and P. Oliveira. 
\newblock Sensor-based globally asymptotically stable filters for attitude estimation: Analysis, design, and performance evaluation. 
\newblock {\em IEEE Transactions on Automatic Control}, 57(8):2095--2100, 2012.
%
\bibitem{bonnabelITAC08}
S. Bonnabel, P. Martin, and P. Rouchon. 
\newblock Symmetry-preserving observers. 
\newblock {\em IEEE Trans. on Autom. Control}, 53(11):2514--2526, 2008.
%
\bibitem{brockett1972}
R. W. Brockett. 
\newblock System theory on group manifolds and coset spaces.
\newblock {\em SIAM Journal on control}, 10(2):265--284, 1972.
%
\bibitem{grip2012nonlinear}
H. F. Grip, T. I. Fossen, T. A. Johansen, and A. Saberi. 
\newblock A nonlinear observer for integration of GNSS and IMU measurements with gyro bias estimation. 
\newblock In {\em American Cont. Conf.}, pages 4607--4612, 2012.
%
\bibitem{HamSamTAC16}
T. Hamel and C. Samson. 
\newblock Riccati observers for the non-stationary PnP problem. 
\newblock Submitted to {\em IEEE Trans. on Automatic Control}, 2016,
available at https://hal.archives-ouvertes.fr/hal-01326999.
%
\bibitem{hua10cep}
M.-D. Hua. 
\newblock Attitude estimation for accelerated vehicles using GPS/INS measurements. 
\newblock {\em Control Eng. Pract.}, 18(7):723--732, 2010.
%
\bibitem{hua14}
M.-D. Hua, G. Ducard, T. Hamel, R. Mahony, and K. Rudin. 
\newblock Implementation of a nonlinear attitude estimator for aerial robotic vehicles.
\newblock {\em IEEE Trans. on Control Syst. Technol.}, 22(1):201--213, 2014.
%
\bibitem{hua16}
M.-D. Hua, P. Martin, and T. Hamel. 
\newblock Stability analysis of velocityaided attitude observers for accelerated vehicles. 
\newblock {\em Automatica}, 63:11--15, 2016.
%
\bibitem{khosravian16}
A. Khosravian, J. Trumpf, R. Mahony, and T. Hamel. 
\newblock State estimation for invariant systems on lie groups with delayed output measurements.
\newblock {\em Automatica}, 68:254--265, 2016.
%
\bibitem{mhp08}
R. Mahony, T. Hamel, and J.-M. Pflimlin. 
\newblock Nonlinear complementary filters on the special orthogonal group. 
\newblock {\em IEEE Trans. on Autom. Control}, 53(5):1203--1218, 2008.
%
\bibitem{mahony2012}
R. Mahony, V. Kumar, and P. Corke. 
\newblock Multirotor aerial vehicles: Modeling, estimation, and control of quadrotor. 
\newblock {\em IEEE Robotics \& Automation Magazine}, 19(3):20--32, 2012.
%
\bibitem{ms08IFAC}
P. Martin and E. Salaun. 
\newblock An invariant observer for Earth-Velocity-Aided attitude heading reference systems. 
\newblock In {\em IFAC World Congr.}, pages 9857--9864, 2008.
%
\bibitem{ms10cep}
P. Martin and E. Salaun. 
\newblock Design and implementation of a low-cost observer-based attitude and heading reference system. 
\newblock {\em Control Eng. Pract.}, 18(7):712--722, 2010.
%
\bibitem{Martin10ICRA}
P. Martin and E. Salaun. 
\newblock The true role of accelerometer feedback in quadrotor control. 
\newblock In {\em IEEE Int. Conf. on Robotics and Automation}, pages 1623--1629, 2010.
%
\bibitem{rt11CDC}
A. Roberts and A. Tayebi. 
\newblock On the attitude estimation of accelerating rigid-bodies using GPS and IMU measurements. 
\newblock In {\em IEEE Conf. on Dec. and Cont.}, pages 8088--8093, 2011.
%
\bibitem{scandaroli2013}
G. G. Scandaroli. 
\newblock Visuo-inertial data fusion for pose estimation and self-calibration. 
\newblock PhD thesis, Universit\'e Nice Sophia Antipolis, 2013.
%
\bibitem{troniICRA13}
G. Troni and L. L. Whitcomb. 
\newblock Preliminary experimental evaluation of a Doppler-aided attitude estimator for improved Doppler navigation of underwater vehicles. 
\newblock In {\em IEEE Int. Conf. on Robotics and Automation}, pages 4134--4140, 2013.
\end{thebibliography}
\end{document}